\documentclass[12pt,draftcls,onecolumn]{IEEEtran}
\usepackage{mathrsfs,amsmath,latexsym,amssymb}

\newtheorem{thm}{Theorem}[section]

\newtheorem{cor}{Corollary}[section]
\newcommand{\argmax}{\mathop{\mbox{\rm arg\,max}}}

\newcommand{\ps}{\mathrm{ps}}
\begin{document}
\sloppy

\title{On Supervised On-line Rolling-Horizon Control for Infinite-Horizon Discounted Markov Decision Processes}
\author{Hyeong Soo Chang
\thanks{H.S. Chang is with the Department of Computer Science and Engineering at Sogang University, Seoul 121-742, Korea. (e-mail:hschang@sogang.ac.kr).}%
}

\maketitle
\begin{abstract}
This note re-visits the rolling-horizon control approach
to the problem of
a Markov decision process (MDP) with infinite-horizon discounted expected reward criterion.
Distinguished from the classical value-iteration approach,
we develop an asynchronous on-line algorithm based on policy iteration integrated with
a multi-policy improvement method of policy switching. 
A sequence of monotonically improving solutions to the forecast-horizon sub-MDP 
is generated by updating the current solution only at the currently visited state, building in effect a rolling-horizon control policy for the MDP 
over infinite horizon.
Feedbacks from ``supervisors," if available, can be also
incorporated while updating.
We focus on the convergence issue with a relation to the transition 
structure of the MDP.
Either a global convergence to an optimal forecast-horizon policy 
or a local convergence to a ``locally-optimal" fixed-policy in a finite time
is achieved by the algorithm depending on the structure.
\end{abstract}


\section{Introduction}

Consider the rolling horizon control (see, e.g.,~\cite{her}) with a fixed finite forecast-horizon $H$ to
the problem of a Markov decision process (MDP) $M_{\infty}$ with infinite-horizon discounted expected
reward criterion.
At discrete time $k\geq 1$, the system is at a state $x_k$ in a finite state set $X$.
If the controller of the system takes an action $a$ in $A(x_k)$ at $k$, then it obtains a reward of $R(x_{k},a)$ from 
a reward function $R:X\times A\rightarrow \Re$, where $A(x)$ denotes an admissible-action set of $x$ in $X$. The system then makes a random transition to a next state $x_{k+1}$ 
by the probability of $P_{x_{k}x_{k+1}}^{a}$ specified in an $|X|\times |X|$-matrix $P^a=[P_{xy}^a], x,y\in X$.

Let $B(X)$ be the set of all possible functions from $X$ to $\Re$. The zero function in $B(X)$ is
referred to as $0_X$ where $0_X(x)=0$ for all $x\in X$.
Let also $\Pi(X)$ be the set of all possible mappings from $X$ to $A$ where for any $\pi \in \Pi(X)$, $\pi(x)$ is constrained to be in $A(x)$ for all $x\in X$.
Let an $h$-\emph{length} \emph{policy} of the controller be an element in $\Pi(X)^h$, $h$-ary Cartesian product, $h\geq 1$. That is, $\pi \in \Pi(X)^h$ is an ordered tuple $(\pi_{1},...,\pi_{h})$ where the $m$th entry of $\pi$ is equal to $\pi_m \in \Pi(X), m \geq 1,$ and 
when $\pi_m$ is applied at $x$ in $X$, then the controller is supposed to look ahead of (or forecasts) the remaining horizon $h-m$ for control.
An infinite-length policy is an element in the infinite Cartesian product of $\Pi(X)$, denoted by 
$\Pi(X)^{\infty}$, and referred to as just a policy.
We say that a policy $\phi \in \Pi(X)^{\infty}$ is \emph{stationary} if $\phi_m=\pi$ for all $m\geq 1$ for some $\pi\in \Pi(X)$. Given $\pi\in \Pi(X)$, $[\pi]$ denotes a stationary policy in $\Pi(X)^{\infty}$ constructed from $\pi$ such that $[\pi]_m=\pi$ for all $m\geq 1$. 

Define the $h$-horizon value function $V^{\pi}_h$ of $\pi$ in $\Pi(X)^H, H \geq h$ such that
\[
V^{\pi}_h(x) = E\left [\sum_{l=1}^{h} \gamma^{l-1} R(X_{l},\pi_{l}(X_l)) \biggl | X_1=x \right ], x\in X, h \in [1,\infty),
\] where $X_l$ is a random variable that denotes a state at the level (of forecast) $l$ by
following the $l$th-entry mapping $\pi_l$ of $\pi$ and a fixed discounting factor $\gamma$ is in $(0,1)$. 

In the sequel, any operator is applied componentwise for the elements in $B(X)$ and in $\Pi(X)$, respectively. Given $\pi\in \Pi(X)^h$ and $x\in X$, $\pi(x)$ is set to be $(\pi_1(x),...,\pi_h(x))$ meaning the ``$x$-coordinate" of $\pi$ here.
As is well known then,
there exists an optimal $h$-length policy $\pi^*(h)$ such that for any $x\in X$,
$V^{\pi^*(h)}_h(x) =\max_{\pi\in \Pi(X)^h} V^{\pi}_h(x)  = V^*_h(x)$, where
$V^*_h$ is referred to as the optimal $h$-horizon value function. 
In particular, $\{V^*_h, h\geq 1\}$ satisfies the optimality principle:
\[
    V^*_h(x) = \max_{a\in A(x)} \biggl ( R(x,a) + \gamma \sum_{y\in X} P_{xy}^a V^{*}_{h-1}(y) \biggr ), x\in X,
\] for any fixed $V^*_0$ in $B(X)$.
Furthermore, $V^*_h$ is equal to the function in $B(X)$ obtained after applying the value iteration (VI) operator $T:B(X)\rightarrow B(X)$ iteratively $h$ times
with the initial function of $V^*_0$: $T(...T(T(T(V^*_0)))) = T^h(V^*_0)) = V^*_h$, where
\[T(u)(x) = \max_{a\in A(x)} ( R(x,a) + \gamma \sum_{y\in X} P_{xy}^a u(y) ), x\in X, u\in B(X).
\]
This optimal substructure leads to a dynamic programming algorithm, backward induction, which computes $\{V^*_h\}$ in \emph{off-line} and returns an optimal $H$-horizon policy $\pi^*(H)$ that achieves the optimal value at any $x\in X$ for the $H$-horizon sub-MDP $M_H$ of $M_{\infty}$ by
\[\pi^*(H)_{H-h+1}(x) \in \argmax_{a\in A(x)} \left ( R(x,a) + \gamma \sum_{y\in X} P_{xy}^a V^{*}_{h-1}(y) \right ), x\in X, h \in \{1,...,H\}.
\]
Once obtained, the rolling $H$-horizon controller employs the first entry $\pi^*(H)_1$ of $\pi^*(H)$ or a stationary policy $[\pi^*(H)_1]$ over the system time: At each $k\geq 1$, $\pi^*(H)_1(x_k)$ is taken at $x_k$.

Given $\pi\in \Pi(X)^{\infty}$,
$V^{\pi}_{\infty}$ denotes the value function of $\pi$ over \emph{infinite horizon} and it is obtained by
letting $h$ approach infinity in $E[\sum_{l=1}^{h} \gamma^{l-1} R(X_{l},\pi_{l}(X_l))|X_1=x]$.
The optimal value function $V^*_{\infty}$ of $M_{\infty}$ is then defined such that $V^*_{\infty}(x) = \sup_{\pi\in \Pi(X)^{\infty}} V^{\pi}_{\infty}(x), x\in X$.
A performance result of $[\pi^*(H)_1]$ applied to $M_{\infty}$ (see, e.g.,~\cite{her}) is that 
the infinity-norm of the difference between the value function of $[\pi^*(H)_1]$ and $V^*_{\infty}$ is upper bounded by (an error of) $O(\gamma^H ||V^*_{\infty} - V^*_0||_{\infty})$. 
The term $||V^*_{\infty} - V^*_0||_{\infty}$ can be loosely upper bounded by $C/(1-\gamma)$ with some constant $C$. Then due to the dependence on $(1-\gamma)^{-1}$, the performance worsens around $\gamma$ closer to one. 
How to set $V^*_0$ is a critical issue in the rolling horizon control
even if the error vanishes to zero exponentially fast in $H$ with the rate of $\gamma$.

This note develops an algorithm for \emph{on-line} rolling $H$-horizon control. The sub-MDP $M_H$ is \emph{not} solved in advance.
Rather with an arbitrarily selected $\pi_1(H)\in \Pi(X)^H$ for $M_H$, the algorithm generates \emph{a monotonically improving sequence of} $\{\pi_{k}(H)\}$ over time $k\geq 1$. To the algorithm, only $\pi_{k-1}(H)$ is available at $k>1$ and it updates
$\pi_{k-1}(H)$ \emph{only at} $x_k$ to obtain $\pi_k(H)$.
Either we have that $\pi_{k}=\pi_{k-1}$ or $\pi_k(x)=\pi_{k-1}(x)$ for all $x\in X\setminus\{x_k\}$ but $\pi_k(x_k)\neq \pi_{k-1}(x_k)$.
The algorithm has a design flexibility in the aspect that a feedback of an action to be used at $x_k$
by some ``supervisor," can be incorporated while generating $\pi_k(H)$.
By setting $\phi_k = \pi_k(H)_1$ at each $k\geq 1$, a policy $\phi \in \Pi(X)^{\infty}$
is in effect built sequentially for the controller.
Once $\phi_k$ is available to the controller, it takes $\phi_{k}(x_k)$ to the system and the underlying system of $M_{\infty}$ moves to a next random state $x_{k+1}$ by the probability of $P^{\phi_k(x_k)}_{x_k x_{k+1}}$.

The behavior of such a control policy is discussed by focusing on the convergence issue with
a relation to
the transition structure of $M_{\infty}$. 
We are concerned with a question about the \emph{existence of a finite time} $K<\infty$ such that
$\phi_{k}=\pi^*(H)_1$ for all $k > K$ for the infinite sequence $\{\phi_k\}$.

\section{Off-line Synchronous Policy Iteration with Policy Switching}
\label{sec:off-line}

To start with, we present an algorithm of \emph{off-line synchronous} policy iteration (PI) combined with 
a multi-policy improvement method of policy switching for solving $M_H$. In what follows, we assume that $V^{\pi}_0 = 0_X$ for any $\pi\in \Pi(X)^H$ for simplicity. (The topic about how to set $V^{\pi}_0$ is beyond the scope of this note.)

\begin{thm}[Theorem 2~\cite{changps}]
Given a nonempty $\Delta \subseteq \Pi(X)^H$, construct policy switching with $\Delta$ in $\Pi(X)^H$ as $\pi_{\ps}(\Delta)$ such that for each possible pair of $x\in X$ and $h \in \{1,...,H\}$,
\[ \pi_{\ps}(\Delta)_h(x) = \phi^*_h(x) \mbox{ where } \phi^* \in \argmax_{\phi\in\Delta} V^{\phi}_{H-h+1}(x).
\] 
Then $V^{\pi_{\ps}(\Delta)}_H\geq V^{\phi}_H$ for all $\phi \in \Delta$.
\end{thm}

Given $\tilde{\pi}$ and $\pi$ in $\Pi(X)^H$,
we say that $\tilde{\pi}$ \emph{strictly improves} $\pi$ (over the horizon $H$)
if there exists some $s\in X$ such that $V^{\tilde{\pi}}_H(s) > V^{\pi}_H(s)$ and $V^{\tilde{\pi}}_H\geq V^{\pi}_H$
in which case we write as $\tilde{\pi} >_H \pi$.
Define \emph{switchable action set} $S^{\pi}_h(x)$ of $\pi\in \Pi(X)^H$ at $x\in X$ for $h\in \{1,...,H\}$ as
\[
S^{\pi}_h(x) = \biggl \{ a \in A(x) \biggl | R(x,a) + \gamma \sum_{y\in X} P^{a}_{xy}V^{\pi}_{h-1}(y) > V^{\pi}_h(x) \biggr \}
\] and also \emph{improvable-state set of} $\pi$ for $h$ as 
\[I^{\pi}_h  = \{ (h,x) | S^{\pi}_h(x) \neq \emptyset, x\in X \}. 
\]
Set $I^{\pi,H} = \bigcup_{h=1}^{H} I^{\pi}_h$.
Observe that if $I^{\pi,H} = \emptyset$ for $\pi\in \Pi(X)^H$, then $\pi$ is an optimal $H$-length policy for $M_H$.

The following theorem provides a result for $M_H$ in analogy with the key step of the single-policy improvement (see, e.g.,~\cite{puterman})
in PI for $M_{\infty}$.
Because Banach's fixed-point theorem is difficult to be invoked in the finite-horizon case unlike the standard proof for the infinite-horizon case, we provide a proof for the completeness.
\begin{thm}
\label{thm:imp}
Given $\pi \in \Pi(X)^H$ with $I^{\pi,H} \neq \emptyset$, construct $\tilde{\pi}\in \Pi(X)^H$ 
with any $I$ satisfying $\emptyset \subsetneq I \subseteq I^{\pi,H}$
such that $\tilde{\pi}_{H-h+1}(x) \in S^{\pi}_h(x)$ for all $(h,x)\in I$ and $\tilde{\pi}_{H-h+1}(x)=\pi_{H-h+1}(x)$ for all $(h,x) \in (\{1,...,H\} \times X)\setminus I$. Then $\tilde{\pi} >_H \pi$.
\end{thm}
\proof
The base case holds because $V^{\tilde{\pi}}_0 = V^{\pi}_0$.
For the induction step, assume that $V^{\tilde{\pi}}_{h-1} \geq V^{\pi}_{h-1}$.
For all $x$ such that $\tilde{\pi}_{H-h+1}(x)=\pi_{H-h+1}(x)$, 
\begin{eqnarray*}
\lefteqn{V^{\tilde{\pi}}_h(x) = R(x,\tilde{\pi}_{H-h+1}(x)) + \gamma \sum_{y\in X} P^{\tilde{\pi}_{H-h+1}(x)}_{xy}V^{\tilde{\pi}}_{h-1}(y)} \\
& & \geq R(x,\pi_{H-h+1}(x)) + \gamma \sum_{y\in X} P^{\pi_{H-h+1}(x)}_{xy}V^{\pi}_{h-1}(y) = V^{\pi}_h(x).
\end{eqnarray*} 
On the other hand, for any $x\in X$ such that $\tilde{\pi}_{H-h+1}(x) \in S^{\pi}_h(x)$,
\begin{eqnarray*}
\lefteqn{V^{\tilde{\pi}}_h(x) = R(x,\tilde{\pi}_{H-h+1}(x)) + \gamma \sum_{y\in X} P^{\tilde{\pi}_{H-h+1}(x)}_{xy}V^{\tilde{\pi}}_{h-1}(y)} \\
& &  \geq R(x,\tilde{\pi}_{H-h+1}(x)) + \gamma \sum_{y\in X} P^{\tilde{\pi}_{H-h+1}(x)}_{xy}V^{\pi}_{h-1}(y)
 \mbox{ by induction hypothesis } V^{\tilde{\pi}}_{h-1} \geq V^{\pi}_{h-1} \\
& & > R(x,\pi_{H-h+1}(x)) + \gamma \sum_{y\in X} P^{\pi_{H-h+1}(x)}_{xy}V^{\pi}_{h-1}(y) = V^{\pi}_h(x)
 \mbox{ because } \tilde{\pi}_{H-h+1}(x)\in S^{\pi}_h(x). 
\end{eqnarray*}  
Putting the two cases together makes $V^{\tilde{\pi}}_{h} \geq V^{\pi}_{h}$. In particular, we see that there must exist some $y\in X$ such that $\tilde{\pi}_{H-h+1}(y) \in S^{\pi}_h(y)$ because of our choice of $I$, having $V^{\tilde{\pi}}_h(y) > V^{\pi}_h(y)$.
By an induction argument on the level from $h$ then, it follows that $V^{\tilde{\pi}}_{H}(y) > V^{\pi}_{H}(y)$.
\endproof

The previous theorem states that if a policy was generated from a given $\pi$ by switching
the action prescribed by $\pi$ at each improvable state with its corresponding level in a
particularly chosen nonempty subset of the improvable-state set of $\pi$, then
the policy constructed strictly improves $\pi$ over the relevant finite horizon.
However, in general, even if $\pi >_H \phi$ is known, for $\pi'$ and $\phi'$ obtained by the method 
of Theorem~\ref{thm:imp}, respectively, $\pi' >_H \phi'$ does \emph{not} hold necessarily.
(Note that this is also true for the infinite-horizon case.)
It can be merely said that $\pi'$ improves $\pi$ and $\phi'$ improves $\phi$, respectively.
Motivated by this, we consider the following. For a given $\pi$ in $\Pi(X)^{H}$, let $\beta^{\pi,H}$ 
be the set of \emph{all strictly better policies than} $\pi$ \emph{obtainable from} $I^{\pi,H}$:
If $I^{\pi,H} = \emptyset$, $\beta^{\pi,H} = \emptyset$. Otherwise,
\begin{eqnarray*}
\lefteqn{\beta^{\pi,H} = \bigl \{\tilde{\pi}\in \Pi(X)^H | I \in 2^{I^{\pi,H}}\setminus \{\emptyset\},
  \forall (h,x)\in I \mbox{ } \tilde{\pi}_{H-h+1}(x) \in S^{\pi}_h(x)} \\ 
& & \hspace{3cm} \mbox{ and } \forall (h,x)\in (\{1,...,H\} \times X) \setminus I \mbox{ } \tilde{\pi}_{H-h+1}(x)=\pi_{H-h+1}(x) \bigr \}.
\end{eqnarray*} Once obtained, policy switching with respect to $\beta^{\pi,H}$ is applied
to find a \emph{single} policy no worse than all policies in the set.

We are ready to derive an off-line synchronous algorithm, ``policy-iteration with policy switching," (PIPS) that generates a sequence of $H$-length policies for solving $M_H$: Set arbitrarily $\pi_1 \in \Pi(X)^H$.
Loop with $n \in \{1,2,...,\}$ until $I^{\pi_{n},H} = \emptyset$ 
where $\pi_{n+1} = \pi_{\ps}(\beta^{\pi_n,H})$.

The convergence to an optimal $H$-length policy for $M_H$ is trivially guaranteed because $\pi_{n+1} >_H \pi_n, n\geq 1$, and both $X$ and $A$ are finite. 
Note that $\beta^{\pi_n,H}$ in $\pi_{\ps}(\beta^{\pi_n,H})$ can be substituted with any $\Delta_n\subseteq \Pi(X)^H$ as long as for $\Delta_n \cap \beta^{\pi_n,H} \neq \emptyset$. 
The generality of $\{\Delta_n\}$ then provides a broad design-flexibility of PIPS.

The idea behind policy switching used in PIPS with $\beta^{\pi_n,H}$ can be attributed to 
approximating the steepest ascent direction 
while applying the steepest ascent algorithm. At the current location $\pi_n$,
we find ascent ``directions" relative to $\pi_n$ over the \emph{local neighborhood of} $\beta^{\pi_n,H}$.
A steepest ascent direction, $\pi_{\ps}(\beta^{\pi_n,H})$, is then obtained by ``combining" all 
of the possible ascent directions. In particular, the \emph{greedy} ascent direction $\phi$ that satisfies that
$T(V^{\pi_n}_{H-h})(x) = R(x,\phi_h(x)) + \gamma \sum_{y\in X} P_{xy}^{\phi_h(x)} V^{\pi_n}_{H-h}(y)$
for all $x\in X$ and $h\in \{1,...,H\}$ is always included while combining.

\section{Off-line Asynchronous PIPS}
\label{sec:on-line}

An \emph{asynchronous} version can be inferred from the synchronous PIPS by the 
following improvement result 
when a single $H$-length policy in $\Pi(X)^H$ is updated only at a \emph{single state}:
\begin{cor}
Given $x\in X$ and $\pi\in \Pi(X)^H$, let $I^{\pi,H}_x = \{(h,x)\in I^{\pi,H} | h\in \{1,...,H\}\}$. 
Suppose that $I^{\pi,H}_x\neq \emptyset$. Then for any $\phi \in \beta^{\pi,H}_{x}$, 
$\phi >_H \pi$ where
\begin{eqnarray*}
\lefteqn{\hspace{-1cm}\beta^{\pi,H}_x = \bigl \{\tilde{\pi}\in \Pi(X)^H | I \in 2^{I^{\pi,H}_x}\setminus \{\emptyset\},
  \forall h \in \{1,...,H\} \mbox{ such that } (h,x)\in I \mbox{ } \tilde{\pi}_{H-h+1}(x)}  \\
  & & \in S^{\pi}_{h}(x) \mbox{ and } \forall (h,x')\in (\{1,...,H\} \times X) \setminus I \mbox{ } \tilde{\pi}_{H-h+1}(x')=\pi_{H-h+1}(x') \bigr \}.
\end{eqnarray*} 
\end{cor}
\proof
Because $\emptyset \neq I^{\pi,H}_x \subseteq I^{\pi,H}$, $\emptyset \neq \beta^{\pi,H}_x \subseteq \beta^{\pi,H}$.
\endproof

The set $\beta^{\pi,H}_x$ ``projected to the $x$-coordinate direction" from $\beta^{\pi,H}$
contains all strictly better policies than $\pi$ 
that can be obtained by switching the action(s) prescribed by $\pi$ at the single state $x$.
This result leads to an off-line convergent asynchronous PIPS for $M_H$:
Select $\pi_1\in \Pi(X)^H$ arbitrarily.
Loop with $n \in \{1,2,...,\}$:
If $I^{\pi_n,H} \neq \emptyset$,
select $x_n \in I^{\pi_n,H}$ and
construct $\pi_{n+1}$ such that $\pi_{n+1}(x_n) = \pi_{\ps}(\beta^{\pi_n,H}_{x_n})(x_n)$ and $\pi_{n+1}(x)=\pi_{n}(x)$ for all $x$ in $X\setminus\{x_n\}$. If $I^{\pi_n,H} = \emptyset$, stop.
Because $x_n$ is always selected to be an improvable-state in $I^{\pi_n,H}\neq \emptyset$, $\pi_{n+1} >_H \pi_n$ for all $n\geq 1$. Therefore, $\{\pi_n\}$ converges to an optimal $H$-length policy for $M_H$.

Suppose that the state given at the current step of the previous algorithm is \emph{not} 
guaranteed to be in the improvable-state set of the current policy.
Such scenario is possible with the following modified version:
Select $\pi_1\in \Pi(X)^H$ arbitrarily.
Loop with $n\in \{1,2,...,\}$:
If $I^{\pi_n,H}=\emptyset$, stop. Select $x_n \in X$. If $x_n \in I^{\pi_n,H}$, then
construct $\pi_{n+1}$ such that $\pi_{n+1}(x_n) = \pi_{\ps}(\beta^{\pi_n,H}_{x_n})(x_n)$ and $\pi_{n+1}(x)=\pi_{n}(x)$ for all $x$ in $X\setminus\{x_n\}$. If $x_n \notin I^{\pi_n,H}$, $\pi_{n+1} = \pi_n$.
Unlike the previous version, this algorithm's convergence
\emph{depends on the sequence $\{x_n\}$ selected}.
Even if $\pi_{n+1} >_H \pi_n$ when $\pi_{n+1}\neq \pi_n$, 
the stopping condition that checks for the optimality can never be satisfied. In other words,
an infinite loop is possible.
The immediate problem is then how to choose an update-state sequence to achieve 
a global convergence.
The reason for bring this issue up with the modified algorithm is that the situation is
closely related with the on-line algorithm to be discussed in the next section. Dealing with this issue
here would help understanding the convergence behaviour of the on-line algorithm.
We discuss some pedagogical example of choosing an update-state sequence of the modified
off-line algorithm below.

One way of \emph{enforcing} a global convergence is to ``embed" backward induction into the update-state sequence. 
For example, we generate a sequence of
$\{x_n\}=\{x_0,...,x_{n_1},...,x_{n_2},...,x_{n_h},....,x_{n_H},...\}$ whose
subsequence $\{x_{n_h},h=1,...,H\}$ produces $\{\pi^{n_h}\}$ that
solves $M_h$.
We need to follow the optimality principle such that $M_{h-1}$ is solved \emph{before} $M_h$, and so forth, until $M_H$ is finally solved. Therefore, the entries of $\pi^*(H)$ are searched from $\pi^*(H)_H$ to $\pi^*(H)_1$ 
over $\{x_n\}$ such that 
$\pi^{n_1} = (\pi^*(H)_H, \pi^{n_1}_2,...,\pi^{n_1}_{H-1},\pi^{n_1}_H)$ where $V^{*}_{1} = V^{\pi^{n_1}}_1$, and then
$\pi^{n_2} = (\pi^*(H)_{H-1}, \pi^*(H)_H,\pi^{n_2}_3,...,\pi^{n_2}_H)$ where $V^{*}_{2} = V^{\pi^{n_2}}_2$,
$\pi^{n_h} = (\pi^*(H)_{H-h+1},...,\pi^*(H)_H,...,\pi^{n_h}_H)$ where $V^{*}_{h} = V^{\pi^{n_h}}_h$,..., and then finally,
$\pi^{n_H} = (\pi^*(H)_1,...,\pi^*(H)_{H-1},\pi^*(H)_H)$ with $V^{*}_{H} = V^{\pi^{n_H}}_H$.

Once $M_{h-1}$ has been solved, an optimal $h$-length policy $\pi^{n_h}$ for $M_h$ can be found exhaustively.
The corresponding update-state subsequence from $x_{n_{h-1}+1}$ to $x_{n_h}$ can be any permutation of the states in $X$. Visiting each $x$ in $X$ \emph{at least once} for updating causes an optimal $h$-length policy for $M_h$ to be found 
because if not empty, $\beta^{\pi^{m},H}_{x}$, where $m\in \{n_{h-1}+1,...,n_h\}$, includes an $H$-length policy
whose first entry mapping maps $x$ to an action in
$\argmax_{a\in A(x)} ( R(x,a) + \gamma \sum_{y\in X} P_{xy}^a V^{*}_{h-1}(y) )$.
Even though visiting each state at least once makes the approach enumerative,
our point is showing that there \emph{exists} an update-state sequence that makes a global convergence 
possible.

\section{On-line Asynchronous PIPS}

We are now at the position of the main subject of this note about solving
$M_{\infty}$ within an on-line rolling $H$-horizon control by solving $M_H$ not in advance but over time.
We assume that a sequence of $\{\Delta_k\}$ is available where $\Delta_k\subseteq \Pi(X)^H, k\geq 1$,
stands for a set of supervisors at $k$.
Any feedback of an action to take at a state can be represented by an element in $\Pi(X)^H$.

The controller applies a policy $\phi\in \Pi(X)^{\infty}$ to the system of $M_{\infty}$
built from the sequence of the $H$-length policies generated by the on-line asynchronous PIPS algorithm: PIPS first selects $\pi_1(H)\in \Pi(X)^H$ arbitrarily
and sets $\phi_1=\pi_1(H)_1$.
At $k >1$, if $\beta^{\pi_{k-1}(H),H}_{x_k} = \emptyset$, $\pi_k(H) = \pi_{k-1}(H)$.
Otherwise, PIPS updates $\pi_{k-1}(H)$ only at $x_k$ such that
$\pi_{k}(H)(x_k) = \pi_{\ps}(\beta^{\pi_{k-1}(H),H}_{x_k} \cup \Delta_k)(x_k)$ 
and $\pi_{k}(H)(x) = \pi_{k-1}(H)(x)$ for all $x$ in $X\setminus \{x_k\}$.
After finishing update, PIPS sets $\pi_{k}(H)_1$ to the $k$th entry $\phi_k$ of $\phi$.
Once $\phi_k$ is available to the controller, $\phi_{k}(x_k)$ is taken to the system and 
the system makes a random transition to $x_{k+1}$.

Before we present a general convergence result, an intuitive consequence about a global convergence
from a sufficient condition related with the transition structure of $M_{\infty}$
is given below as a theorem.
Note that $M_{\infty}$ is communicating, if every Markov chain 
induced by fixing each stationary policy $[\pi], \pi \in \Pi(X),$ in $M_{\infty}$ is communicating~\cite{kallen}.
\begin{thm}
\label{thm:mainconv}
Suppose that $M_{\infty}$ is communicating. Then for $\{\pi_k(H), k\geq 1\}$ generated 
by on-line asynchronous PIPS,
there exists some $K < \infty$ such that $\pi_k(H) = \pi^*(H)$ for all $k > K$ where 
$V^{\pi^*(H)}_H = V^*_H$ in $M_H$
for any $\pi_1(H)\in \Pi(X)^H$ and $\{\Delta_k\}$ where $\Delta_k\subseteq \Pi(X)^H, k\geq 1$.
\end{thm}
\proof
Because both $X$ and $A$ are finite, $B(X)$ and $\Pi(X)^H$ are finite. For any given $\{\Delta_k\}$
and any $\pi_1(H)$,
the monotonicity of $\{V^{\pi_k(H)}\}$ of $\{\pi_k(H)\}$ holds because
$\{\pi_k(H)\}$ satisfies that $\pi_k(H) >_H \pi_{k-1}(H)$ if $\beta^{\pi_{k-1}(H),H}_{x_k} \neq \emptyset$ and $V^{\pi_k(H)} \geq V^{\pi_{k-1}(H)}$ otherwise.
The assumption that $M_{\infty}$ is communicating ensures that every state $x$ in $X$ is visited 
infinitely often within $\{x_k\}$. It follows that at some sufficiently large finite time $K$, 
$\beta^{\pi_{K}(H),H} = \emptyset$ implying $I^{\pi_{K}(H),H}=\emptyset$. Therefore, $\pi_k(H) = \pi^*(H)$ for all $k > K$ where $V^{\pi^*(H)}_H = V^*_H$ in $M_H$.
\endproof
The policy $\phi$ of the controller becomes \emph{stable}
in the sense that the sequence $\{\phi_k\}$ converges to $\pi^*(H)_1$. 
The question about checking whether $M_{\infty}$ is communicating without enumerating all stationary
policies in $\Pi(X)^{\infty}$ is possible in a polynomial time-complexity was raised
in~\cite{kallen}. Unfortunately, this problem is in general NP-complete~\cite{tsi}. 
A simple and obvious sufficient condition for such connectivity is that $P^{a}_{xy} > 0$ for 
all $x,y$ in $X$ and $a$ in $A(x)$.
The key in the convergence here is that which state in $X$ is visited ``sufficiently often" by 
following $\{\phi_k\}$ to ensure that an optimal action at the visited state is eventually found.
The following result stated in Theorem~\ref{thm:gen} reflects this rationale.
Given a stationary policy $\phi\in \Pi(X)^{\infty}$, 
the \emph{connectivity relation} $\chi^{\phi}$ is defined on $X$ from the Markov chain $M_{\infty}^{\phi}$
induced by fixing $\phi$ in $M_{\infty}$: If $x$ and $y$ in $X$ communicate each other in $M_{\infty}^{\phi}$,
$(x,y)$ is an element of $\chi^{\phi}$.
Given $x\in X$, the equivalence class of $x$ 
with respect to $\chi^{\phi}$ is denoted by $[x]_{\chi^{\phi}}$. Note that
for any $x\neq y$, either $[x]_{\chi^{\phi}} = [y]_{\chi^{\phi}}$ or  
$[x]_{\chi^{\phi}} \cap [y]_{\chi^{\phi}} = \emptyset$. The collection
of $[x]_{\chi^{\phi}}, x\in X,$ partitions $X$.
\begin{thm}
\label{thm:gen}
For any $\pi_1(H)\in \Pi(X)^H$ and any $\{\Delta_k\}$ where $\Delta_k\subseteq \Pi(X)^H, k\geq 1$,
$\{\pi_k(H)\}$ generated by on-line asynchronous PIPS converges to some $\lambda(H)$ in $\Pi(X)^H$
such that for some $K < \infty$, $\pi_k(H) = \lambda(H)$ for all $k > K$.
Furthermore, $\lambda(H)$ satisfies that 
$V^{\lambda(H)}_H \geq V^{\pi}_H$
for all $\pi \in \bigcup_{x\in [x^*]_{\chi^{[\lambda(H)_1]}} } \beta^{[\lambda(H)_1],H}_x$, where $x^*$
is any visited state at $k>K$.
\end{thm} 
\proof
By the same reasoning in the proof of Theorem~\ref{thm:mainconv}, $\{\pi_k(H)\}$
converges to an element $\lambda(H)$ in $\Pi(X)^H$ in a finite time $K$. 
Because every state $x$ in $[x^*]_{\chi^{[\lambda(H)_1]}}$ is visited infinitely 
often within $\{x_k\}$ for $k > K$, $I^{\lambda(H),H}_x = \emptyset$ for such $x$. No more improvement 
is possible at all states in $[x^*]_{\chi^{[\lambda(H)_1]}}$. Otherwise, it contradicts
the convergence to $\lambda(H)$.
\endproof
The theorem states that $\lambda(H)$ is ``locally optimal" over $[x^*]_{\chi^{[\lambda(H)_1]}}$
in the sense that no more improvement is possible at all states in $[x^*]_{\chi^{[\lambda(H)_1]}}$.

We remark that the above local convergence result is different from that of Proposition 2.2 
by Bertsekas~\cite{bert2}. 
In our case, a subset of $X$ in which every state is visited infinitely often 
is \emph{not} assumed to be given in advance.
The sequence of policies generated by on-line PIPS 
will eventually converge to a policy and Theorem~\ref{thm:gen}
characterizes its local optimality with respect to the communicating classes,
in which every state is visited infinitely often, induced by the policy.
Note further that Bertsekas' result is within the context of rolling an \emph{infinite-horizon} 
control.

Unfortunately, we cannot provide a useful performance result about the performance of $\lambda(H)$ 
in Theorem~\ref{thm:gen} here, e.g, an upper bound on $||V^{[\lambda(H)_1]}_{\infty} - V^*_{\infty}||_{\infty}$ 
because it is difficult to bound $||V^{\lambda(H)}_{H} - V^*_{H}||_{\infty}$.
The degree of approximation by $\lambda(H)$ for $\pi^*(H)$ will determine
the performance of the rolling horizon control policy by the on-line asynchronous
PIPS algorithm.

\section{Concluding Remarks}

The off-line and on-line PIPS algorithms can play the role of \emph{frameworks} 
for solving MDPs with supervisors. 
While the disposition of the algorithms (and their developments) was done 
mainly in the perspective of theoretical soundness and results,
both can be easily implemented by simulation if the MDP model is \emph{generative}.
In particular, for the on-line case, each $\pi \in \beta^{\pi_k,H}_{x_k}\cup \Delta_k$
is simply followed (rolled out) over a relevant forecast-horizon 
(see, e.g.,~\cite{bert1} for related approaches and references) in order to
generate its sample-paths. 
If $\beta^{\pi_k,H}_{x_k}\cup \Delta_k$ is large, some policies from $\beta^{\pi_k,H}_{x_k}$
and $\Delta_k$ can be ``randomly" sampled, for example, without losing the improvement 
of $\pi_k$.
A study on the actual implementation is important and is left as a future study.

On-line PIPS is also in the category of ``learning" control.
Essentially, $V^*_0$ can be thought as an initial knowledge of control to the 
system, e.g., as in the value function of a self-learned Go-playing policy 
of AlphaZero from a neural-network based learning-system.
The monotonically improving value-function sequence generated by PIPS 
can be interpreted as learning or obtaining a better knowledge about
controlling the system,

There exist an algorithm, ``adaptive multi-stage sampling," 
(AMS)~\cite{changams} for $M_H$ whose random estimate converges to $V^*_H(x)$
for a given $x$
as the number of samplings approach infinity in the expected sense.
AMS employed within the rolling-horizon control is closer to its original spirit,
compared with on-line PIPS, because the value of $V^*_H(x_k)$ is approximated 
at each visited $x_k$ like solving $M_H$ in advance.
In contrast with the PI-based approach here, the idea of AMS is to replace the 
inside of the maximization over the admissible action set in the $T$-operator 
of VI such that the maximum selection is done with a different utility 
for each action or the ``necessity" measure of sampling, which
is estimated over a set of currently sampled next states with a support
from a certain upper-confidence-bound that controls a degree of search by
the action.
Because the replacement is applied at each level while emulating the process
of backward induction, AMS requires a recursive process in a depth-first-search manner
that effectively builds a sampled tree whose size is exponential in $H$.
It is therefore not surprising that similar to the complexity of backward induction,
AMS also has the (sample) complexity that \emph{exponentially} depends on $H$.
On the other hand, in general estimating the value of a policy, \emph{not the optimal value}, 
is a much easier task by simulation.
Generating random sample-paths of a policy has a \emph{polynomial} dependence on $H$.
More importantly, it seems difficult to discuss the convergence behaviour  
of the rolling horizon AMS-control 
because some characterization of a \emph{stochastic} policy, due to random 
estimates of $V^*_H(x_k)$ at each $k$, needs to be made with a finite sampling-complexity.
Arguably, this would be true for any algorithm that produces a \emph{random} estimate of 
the optimal value, e.g., Monte-Carlo Tree Search (MCTS) used in AlphaGo or 
AlphaZero~\cite{bert1}. However, it is worthwhile to note that these algorithms' output can
act as a feedback from a supervisor in the framework of on-line PIPS.

It can be checked that another multi-policy improvement method 
of parallel rollout~\cite{changps} does not work
for preserving the monotonicity property with \emph{asynchronous} update
when the set of \emph{more than one} policies is applied to the method for the improvement.
Even with synchronous update, the parallel-rollout approach 
requires estimating a \emph{double expectation} for each action, one for the next-state 
distribution and another one in the value function (to be evaluated at the next state).
In contrast, in policy switching a single expectation for each policy needs to be estimated 
leading to a lower simulation-complexity.

\end{document}